\documentclass[a4paper]{article}
\usepackage{amsmath}
\usepackage{amscd}
\usepackage{amstext}
\usepackage{amsfonts}
\usepackage{euscript}
\usepackage{graphicx}

\def\scr{\EuScript}

\newcommand{\pcirc}{{\scriptstyle \,\circ\,}}

\newcommand{\C}{\mathbb{C}}
\newcommand{\ZZ}{\mathbb{Z}}

\DeclareMathOperator{\Hom}{\it Hom}

\DeclareMathOperator{\End}{\it End}

\DeclareMathOperator{\Ext}{\it Ext}

\DeclareMathOperator{\divi}{\rm div}

\DeclareMathOperator{\Id}{\rm Id}

\DeclareMathOperator{\Der}{\it Der}

\DeclareMathOperator{\Tens}{{\rm T}^\bullet_k}
\DeclareMathOperator{\Sim}{{\rm Sym}^{\bullet}_{\OX}}
\DeclareMathOperator{\U}{{\rm U}}

\newcommand{\derlogD}{\Der(\log D)}
\newcommand{\VCERO}{{\DX(\log D)}}
\newcommand{\VO}{{\cal V}_0}
\newcommand{\OmX}{\Omega_X^{\bullet}}

\newcommand{\D}{{\scr D}}
\newcommand{\K}{{\scr K}}
\newcommand{\E}{{\scr E}}
\renewcommand{\L}{{\scr L}}
\newcommand{\M}{{\scr M}}
\newcommand{\N}{{\scr N}}
\newcommand{\PP}{{\scr P}}
\newcommand{\QQ}{{\scr Q}}
\newcommand{\OO}{{\scr O}}

\DeclareMathOperator{\Dual}{\Bbb D}
\newcommand{\Lotimes}{\stackrel{L}{\otimes}}
\DeclareMathOperator{\DR}{DR}
\newcommand{\cd}[3]{{D}^{#1}_{#2}(#3)}

\DeclareMathOperator{\Gr}{Gr}

\newcommand{\OX}{{\scr O}_X}
\newcommand{\DX}{{\scr D}_X}

\newcommand{\SP}{Sp^{\bullet}}
\newcommand{\SPlog}{\SP_{\VCERO}}

\newcommand{\dx}{{\partial_x}}
\newcommand{\dz}{{\partial_z}}
\newcommand{\dy}{{\partial_y}}

\newcounter{numero}[section]
\renewcommand{\thenumero}{(\thesection .\arabic{numero})}
\newcounter{nnumero}[numero]
\renewcommand{\thennumero}{(\thesection .\arabic{numero}.\arabic{nnumero})}
\newenvironment{corolario}{\medskip
\refstepcounter{numero}\noindent {\sc  \thenumero\ Corollaire.}\
\it}{\vspace{1ex}\par}

\newenvironment{ncorolario}{\medskip
\refstepcounter{nnumero}\noindent {\sc  \thennumero\ Corollaire.}\
\it}{\vspace{1ex}\par}

\newenvironment{teorema}{\medskip
\refstepcounter{numero}\noindent {\sc  \thenumero\ Th\'eor\`eme.}\
\it}{\vspace{1ex}\par}

\newenvironment{nteorema}{\medskip
\refstepcounter{nnumero}\noindent {\sc  \thennumero\ Th\'eor\`eme.}\
\it}{\vspace{1ex}\par}

\newenvironment{lema}{\medskip
\refstepcounter{numero}\noindent {\sc  \thenumero\ Lemme.}\
\it}{\vspace{1ex}\par}

\newenvironment{nlema}{\medskip
\refstepcounter{nnumero}\noindent {\sc  \thennumero\ Lemme.}\
\it}{\vspace{1ex}\par}

\newenvironment{ndefinicion}{\medskip
\refstepcounter{nnumero}\noindent {\sc  \thennumero\ D\'efinition.}\
\it}{\vspace{1ex}\par}

\newenvironment{proposicion}{\medskip
\refstepcounter{numero}\noindent {\sc  \thenumero\ Proposition.}\
\it}{\vspace{1ex}\par}

\newenvironment{nproposicion}{\medskip
\refstepcounter{nnumero}\noindent {\sc  \thennumero\ Proposition.}\
\it}{\vspace{1ex}\par}

\newenvironment{nota}{\medskip
\refstepcounter{numero}\noindent {\sc  \thenumero\ Remarque.}\
}{\vspace{1ex}\par}

\newenvironment{nnota}{\medskip
\refstepcounter{nnumero}\noindent {\sc  \thennumero\ Remarque.}\
}{\vspace{1ex}\par}

\newenvironment{ejemplo}{\medskip
\refstepcounter{numero}\noindent {\sc  \thenumero\ Exemple.}\
}{\vspace{1ex}\par}

\newenvironment{problema}{\medskip
\refstepcounter{numero}\noindent {\sc  \thenumero\ Probl\`eme.}\
}{\vspace{1ex}\par}

\newcommand\numero{\medskip\refstepcounter{numero}\noindent{\sc \thenumero}\hspace{1em}}

\newenvironment{prueba}{
\noindent {\sc  Preuve.}\ }{\hfill $\Box$\vspace{1ex}\par}

\title{Dualit\'e et comparaison sur les complexes de de Rham logarithmiques par rapport aux diviseurs
libres}

\author{F. J. Calder\'on Moreno et L. Narv\'aez Macarro\thanks{
The authors are partially supported by BFM2001-3207 and FEDER.}}
\date{}

\begin{document}
\maketitle

\section*{Introduction}

Soit $X$ une vari\'et\'e analytique complexe lisse de dimension $n$
et $D\subset X$ un diviseur ( $=$ hypersurface) libre, i.e. un
diviseur tel que le faisceau $\derlogD$  des champs de vecteurs
logarithmiques par rapport \`a $D$ est un $\OX$-module localement
libre (de rang $n$) \cite{ksaito_log}.

Les diviseurs lisses, les diviseurs \`a croisements normaux, les
courbes planes, la r\'eunion des hyperplans de r\'eflexion d'un
groupe de r\'eflexions complexes, les discriminants des applications
stables et les vari\'et\'es de bifurcation sont des exemples de
diviseurs libres.

Comme dans le cas des croisements normaux \cite{del_70}, on d\'efinit
une connexion logarithmique par rapport \`a $D$ comme un $\OX$-module
localement libre $\E$ muni d'une connexion \`a coefficients dans les
1-formes logarithmiques
$$ \nabla': \E\xrightarrow{} \E\otimes_{\OX}\Omega^1_X(\log D),$$
ou ce qui revient au m\^eme, d'un morphisme $\OX$-lin\'eaire \`a
gauche
$$\nabla:\derlogD \xrightarrow{} \End_{\C_X}(\E)$$qui satisfait la r\`egle de Leibniz
$\nabla(\delta)(ae) = a\nabla(\delta)(e) + \delta(a)e$.
L'int\'egrabilit\'e de $(\E,\nabla')$ est caract\'eris\'ee par le
fait que $\nabla$ respecte le crochet de Lie.

Toute connexion logarithmique (par rapport \`a $D$) int\'egrable a un
complexe de de Rham associ\'e $\Omega^{\bullet}_X(\log D)(\E)$.

Dans \cite{es_vi_86}, et toujours dans le cas des croisements
normaux, on \'etudie les connexions logarithmiques int\'egrables
comme modules sur un certain faisceau d'anneaux d'op\'erateurs
diff\'erentiels. En particulier, on interpr\`ete le complexe de de
Rham logarithmique d'une connexion logarithmique int\'egrable par ce
moyen et on d\'ecrit son dual de Verdier.

Dans \cite{calde_ens} le premier auteur a g\'en\'eralis\'e le point
de vue de \cite{es_vi_86} au cas des diviseurs libres arbitraires. En
particulier il a d\'emontr\'e que le terme $0$ de la $V$-filtration
de Malgrange-Kashiwara par rapport \`a $D$ sur le faisceau
d'op\'erateurs diff\'erentiels $\DX$, not\'e $\VCERO$, est
l'alg\`ebre enveloppante de l'alg\'ebro\"{\i}de de Lie $\Der(\log D)$ dans
le sens de \cite{rine-63}. Son gradu\'e pour la filtration par
l'ordre s'identifie \`a l'alg\`ebre sym\'etrique du $\OX$-module
localement libre $\Der(\log D)$. Donc $\VCERO$ est un faisceau
d'anneaux coh\'erent \`a fibres noeth\'eriennes de dimension
homologique finie (voir aussi \cite{backe_96}). En particulier, se
donner une structure de $\VCERO$-module \`a gauche (resp. \`a droite)
sur un $\OX$-module $\M$ est \'equivalent \`a se donner un morphisme
$\OX$-lin\'eaire $\nabla:\Der(\log D) \to \End_{\C_X}(\M)$ tel que
$\nabla(a\delta)(m) = a\nabla(\delta)(m)$, $\nabla([\delta,\delta'])
= [\nabla(\delta),\nabla(\delta')]$ et qui satisfait la r\`egle de
Leibniz $\nabla(\delta)(am) = a\nabla(\delta)(m) + \delta(a)m$ (resp.
est \'equivalent \`a se donner un morphisme $\OX$-lin\'eaire $\nabla
:\Der(\log D) \to \End_{\C_X}(\M)$ tel que $\nabla(a\delta)(m) =
\nabla(\delta)(am)$, $\nabla([\delta,\delta']) =
-[\nabla(\delta),\nabla(\delta')]$ et $\nabla(\delta)(am) =
a\nabla(\delta)(m) - \delta(a)m$), pour toutes les sections locales
$a$ de $\OX$, $\delta,\delta'$ de $\Der(\log D) $ et $m$ de $\M$.

Dans cet article on donne un th\'eor\`eme de dualit\'e qui \'echange
la dualit\'e dans le sens des $\DX$-modules avec une dualit\'e
``tordue" dans le sens des $\VCERO$-modules pour les connexions
logarithmiques int\'egrables par rapport \`a un diviseur libre
arbitraire (voir th. \ref{main} et cor. \ref{cor-main}). L'existence
d'un tel th\'eor\`eme a \'et\'e motiv\'ee par les travaux
\cite{es_vi_86} (appendice) et \cite{ucha_tesis,cas_ucha_stek}.
 D'un point de vue technique, le
r\'esultat pr\'ec\'edent utilise une formule d'associativit\'e pour
des produits tensoriels mixtes sur deux alg\`ebres enveloppantes
embo\^{\i}t\'ees que nous n'avons pas trouv\'e dans la litt\'erature
(voir th. \ref{teor-fund} et cor. \ref{A.2}).

Comme application du th\'eor\`eme pr\'ec\'edent nous d\'ecrivons le
dual de Verdier du complexe de de Rham logarithmique d'une connexion
logarithmique int\'egrable dans le cas o\`u le diviseur $D$ est de
Koszul (\cite{calde_ens}, def. 4.1.1), ce qui g\'en\'eralise le
r\'esultat d\'ej\`a mentionn\'e de \cite{es_vi_86} pour les diviseurs
\`a croisements normaux. Toutes les courbes planes et tous les
diviseurs libres localement quasi-homog\`enes sont de Koszul
\cite{calde_ens,calde_nar_LQHKF}. Nous donnons aussi un crit\`ere
pour la perversit\'e des complexes de de Rham logarithmiques.

Ensuite, nous donnons une caract\'erisation diff\'erentielle du
th\'eor\`eme de comparaison logarithmique (TCL) de
\cite{calde_mond_nar_cas} (voir cor. \ref{criterio}), ainsi qu'une
preuve diff\'erentielle de ce th\'eor\`eme dans le cas des diviseurs
libres localement quasi-homog\`enes, qui avait \'et\'e d\'emontr\'e
dans \cite{cas_mond_nar_96} par une voie topologique. Cette preuve
s'appuie dans des arguments de \cite{cas_ucha_stek} et r\'epond \`a
la conjecture \'enonc\'ee \`a la page 94 de {\it loc.~cit.}. Des
r\'esultats proches ont \'et\'e obtenus dans \cite{torre-45-bis} pour
le cas des diviseurs de Koszul.

Une bonne partie de nos r\'esultats s'\'etendent sans peine au cas
g\'en\'eral des alg\'ebro\"{\i}des de Lie (cf.
\cite{rine-63,MR99b:17021,MR2000f:53109,MR2000m:32015}), mais nous
avons pr\'ef\'er\'e de rester dans le cadre logarithmique, o\`u se
trouve notre motivation originale.

Voici le contenu de cet article. Dans la section 1 on d\'efinit les
complexes de Spencer logarithmiques pour les connexions logaritmiques
int\'egrables.

Dans la section 2 on rappelle tout d'abord des constructions bien
connues sur la structure de $\VCERO$-modules sur le $\otimes_{\OX}$
et le $\Hom_{\OX}$ des $\VCERO$-modules, et on d\'emontre le
th\'eor\`eme \ref{teor-fund} qui joue un r\^ole fondamental dans les
r\'esultats de la section suivante.

Dans la section 3 on d\'emontre les r\'esultats principaux de cet
article: la formule qui \'echange les dualit\'es au niveau de
$\VCERO$ et de $\DX$ pour les connexions logarithmiques
int\'egrables, et la formule qui d\'ecrit le dual de Verdier du
complexe de de Rham logarithmique d'une connexion logarithmique
int\'egrable par rapport \`a un diviseur Koszul-libre.

La section 4  concerne le TCL. On donne un crit\`ere diff\'erentiel
pour le TCL valable pour tout diviseur libre. Ensuite on donne une
nouvelle preuve de nature diff\'erentielle du TCL pour les diviseurs
localement quasi-homog\`enes.

Dans la section 5 on \'etudie un exemple de diviseur libre en
dimension 3, trait\'e dans \cite{cas_ucha_pams}, et on montre que son
complexe de de Rham logarithmique n'est pas pervers, ce qui r\'epond
\`a une question pos\'ee dans \cite{calde_nar_compo}. Finalement on
propose quelques probl\`emes.

Dans l'appendice on se place dans le cadre g\'en\'eral des alg\`ebres
enveloppantes des $(k,A)$-alg\`ebres de Lie et on donne les
r\'esultats n\'ecessaires pour d\'emontrer le th\'eor\`eme
\ref{teor-fund}.

Nous remercions F. J. Castro Jim\'enez, T. Torrelli et J. M. Ucha qui
nous ont expliqu\'e leurs r\'esultats. Nous remercions aussi M.
Schulze par ses remarques sur une version pr\'eliminaire de cet
article.

\section{Modules sur le faisceau des op\'erateurs diff\'e\-rentiels logarithmiques}
\label{sec:1}

 \numero \label{nume-SP} {\em Complexe de
Cartan-Eilenberg-Chevalley-Rinehart-Spencer logarithmique}

 Le faisceau structural $\OX$ est un $\DX$-module \`a
gauche qui, par restriction de scalaires, est aussi un
$\VCERO$-module \`a gauche coh\'erent. En fait on dispose d'une
r\'esolution du type Cartan-Eilenberg-Chevalley-Rinehart-Spencer
(\cite{rine-63}, 4; \cite{calde_ens}, 3.1), not\'ee $\SP_{\VCERO}$ et
d\'efinie de la fa\c{c}on suivante:
$$Sp^{-k}_{\VCERO} = \VCERO\otimes_{\OX}
\stackrel{k}{\bigwedge}\derlogD,\quad k=0,\dots, n$$ et la
diff\'erentielle $d^{-k}: Sp^{-k}_{\VCERO} \to Sp^{-k+1}_{\VCERO}$
est donn\'ee par:
\begin{eqnarray*}
&d^{{-1}}(P\otimes\delta ) = P\delta, &\\
&d^{{-k}}( P\otimes(\delta_1\wedge\cdots\wedge\delta_k))
 =\sum_{i=1}^k
(-1)^{i-1} P\delta_i\otimes(\delta_1\wedge\cdots\widehat{\delta_i}
\cdots\wedge\delta_k)&\\ &+ \sum_{1\leq i<j\leq k}
(-1)^{i+j}P\otimes([\delta_i,\delta_j]\wedge\delta_1\wedge\cdots
\widehat{\delta_i}\cdots\widehat{\delta_j}
       \cdots\wedge\delta_k), \  2\leq k\leq n.&
\end{eqnarray*}
Le morphisme
$$ P\in \VCERO = Sp^{0}_{\VCERO}\mapsto d^0(P):= P(1)\in\OX$$
fait de $\SP_{\VCERO}$ une r\'esolution localement libre de $\OX$
comme $\VCERO$-module. Pour le voir on proc\`ede comme dans {\it
loc.~cit.} (voir \cite{calde_ens}, th. 3.1.2): on filtre le complexe
augment\'e $\SP_{\VCERO}\xrightarrow{d^0} \OX$ par
$$ F^i Sp^{-k}_{\VCERO} = (F^{i-k}\VCERO) \otimes_{\OX}
\stackrel{k}{\bigwedge}\derlogD,\quad F^i\OX = \OX$$ pour $i\geq 0$
et $k=0,\dots,n$, de mani\`ere que le gradu\'e associ\'e est
canoniquement isomorphe au complexe augment\'e
$$\Sim (\derlogD) \otimes_{\OX}
\bigwedge^{\bullet} \derlogD \xrightarrow{d^0} \OX,$$ dont la
diff\'erentielle est donn\'ee par
$$ d^{-k} (F\otimes(\delta_1\wedge\cdots\wedge\delta_k)) =\sum_{i=1}^k (-1)^{i-1}
(F\delta_i)\otimes(\delta_1\wedge\cdots\wedge\widehat{\delta_i}
\wedge\cdots\wedge\delta_k)$$pour $F\in \Sim (\derlogD)$,
$\delta_1,\dots,\delta_k\in \derlogD$ et $k=1,\dots, n$, et
l'augmentation $$d^0: \Sim(\derlogD) \otimes_{\OX}
\stackrel{0}{\bigwedge} \derlogD \to \OX$$ provient de la structure
naturelle de $\OX$ comme module sur l'alg\`ebre sym\'etrique
$\Sim(\derlogD)$. Or, ce complexe augment\'e $$\Sim (\derlogD)
\otimes_{\OX} \bigwedge^{\bullet} \derlogD \xrightarrow{d^0} \OX$$
est exact (cf. \cite{bou_a_10}, \S 9, 3) et on d\'eduit le r\'esultat
cherch\'e. \medskip

\numero {\it Connexions logarithmiques int\'egrables}
\smallskip

 Le faisceau des fonctions m\'eromorphes \`a p\^oles le long de $D$, $\OX(\star D)$, est un
$\DX$-module \`a gauche et donc un $\VCERO$-module \`a gauche.

Pour chaque entier $m$, notons par $\OX(mD)$ le sous-$\OX$-module
localement libre de rang 1 de $\OX(\star D)$ form\'e par les
fonctions m\'eromorphes $h$ telles que $\divi (h) + mD \geq 0$. Il
est clair que chaque $\OX(mD)$ est un $\VCERO$-module \`a gauche.

Si $f=0$ est une \'equation locale r\'eduite de $D$ au voisinage d'un
point $p\in D$, alors $f^{-m}$ est une base locale de $\OO_{X,p}(mD)$
comme $\OO_{X,p}$-module. De plus, si $\delta_1,\dots,\delta_n$ est
une base locale de $\derlogD_p$ et $\delta_i(f)=\alpha_i f$, alors en
vertu de \cite{calde_ens}, th. 2.1.4, on a une pr\'esentation locale
\begin{equation}\label{eq:presenta}
\OO_{X,p}(mD) \simeq \D_{X,p}(\log D)/\D_{X,p}(\log
D)(\delta_1+m\alpha_1,\dots, \delta_n+m\alpha_n).\end{equation}

\begin{ndefinicion} Une connexion logarithmique int\'egrable (le long de
$D$) est un $\VCERO$-module \`a gauche qui est localement libre de
rang fini comme $\OX$-module.
\end{ndefinicion}

Les $\OX(mD)$ sont des connexions logarithmiques int\'egrables.
\smallskip

Toute connexion logarithmique int\'egrable $\E$ est un
$\VCERO$-module coh\'erent. En fait, on peut exhiber comme dans
\ref{nume-SP} une r\'esolution de
Cartan-Eilenberg-Chevalley-Rinehart-Spencer $\SP_{\VCERO}(\E)$ avec
$$Sp^{-k}_{\VCERO}(\E) = \VCERO\otimes_{\OX}
\stackrel{k}{\bigwedge}\derlogD\otimes_{\OX} \E,\quad k=0,\dots, n$$
et la diff\'erentielle $\varepsilon^{-k}: Sp^{-k}_{\VCERO}(\E) \to
Sp^{-k+1}_{\VCERO}(\E)$ donn\'ee par:
\begin{eqnarray*}\displaystyle
&\varepsilon^{-1}(P\otimes\delta\otimes e) = P\delta\otimes e - P\otimes \delta e, &\\
&\varepsilon^{-k}(
P\otimes(\delta_1\wedge\cdots\wedge\delta_k)\otimes e)
 =\sum_{i=1}^k
(-1)^{i-1} P\delta_i\otimes(\delta_1\wedge\cdots\widehat{\delta_i}
\cdots\wedge\delta_k)\otimes e&\\ \displaystyle& - \sum_{i=1}^k
(-1)^{i-1} P\otimes(\delta_1\wedge\cdots\widehat{\delta_i}
\cdots\wedge\delta_k)\otimes (\delta_i e)&\\ \displaystyle &+
\sum_{1\leq i<j\leq k}
(-1)^{i+j}P\otimes([\delta_i,\delta_j]\wedge\delta_1\wedge\cdots
\widehat{\delta_i}\cdots\widehat{\delta_j}
       \cdots\wedge\delta_k)\otimes e, \  2\leq k\leq n.&
\end{eqnarray*}
Le complexe $\SP_{\VCERO}(\E)$ est augment\'e par
$$ P\otimes e\in \VCERO \otimes_{\OX}\E \mapsto \varepsilon^0(P\otimes e):= Pe\in\E.$$

Pour montrer que $\SP_{\VCERO}(\E)$ est une r\'esolution  de $\E$
comme $\VCERO$-module on proc\`ede de mani\`ere tout \`a fait
analogue \`a \ref{nume-SP} en consid\'erant toujours la filtration
constante sur $\E$. Le gradu\'e associ\'e du complexe augment\'e
$\SP_{\VCERO}(\E) \xrightarrow{\varepsilon^0} \E$ est alors
canoniquement isomorphe au tensoris\'e par $\E$ sur $\OX$ du gradu\'e
associ\'e du complexe augment\'e $\SP_{\VCERO}
\xrightarrow{\varepsilon^0} \OX$, et donc il est exact.
\medskip

D'apr\`es \cite{calde_ens}, th. 3.2.1, nous savons que pour tout
$\VCERO$-module \`a gauche $\E$, on a un isomorphisme naturel
\begin{equation} \label{eq:calde-ens-3.2.1}
\Omega^{\bullet}_X(\log D)(\E)\simeq \Hom_{\VCERO}(\SP_{\VCERO},\E)
\simeq R\Hom_{\VCERO}(\OX,\E).\end{equation}

\begin{ndefinicion} \label{def-admi} Nous dirons qu'une connexion logarithmique int\'egrable  $\E$ est admissible si
le complexe $\DX\Lotimes_{\VCERO} \E$ est concentr\'e en degr\'e $0$
et $\DX\otimes_{\VCERO} \E$ est un $\DX$-module holonome.
\end{ndefinicion}

Dire que $\OX$ est une connexion logarithmique int\'egrable
admissible revient \`a dire que $D$ est un {\em diviseur de Spencer}
dans la terminologie de \cite{cas_ucha_stek}, d\'ef. 3.3.

La proposition suivante est une g\'en\'eralisation de
\cite{calde_ens} prop. 4.1.3, th. 4.2.1. et se d\'emontre de
mani\`ere tout \`a fait analogue \`a celle-ci en consid\'erant le
complexe $\SP_{\VCERO}(\E)$:

\begin{nproposicion} \label{K-admi} Supposons que $D$ est un diviseur
libre de Koszul (\cite{calde_ens}, d\'ef. 4.1.1). Alors, toute
connexion logarithmique int\'egrable  est admissible.
\end{nproposicion}

Notons $\omega_X$ (resp. $\omega_X(\log D)$) le faisceau des
$n$-formes diff\'erentielles holomorphes (resp. \`a p\^oles
logarithmiques le long de $D$) sur $X$. Il s'agit du $\OX$-dual du
module $\stackrel{n}{\bigwedge} \Der_{\C_X}(\OX)$ (resp. du module
$\stackrel{n}{\bigwedge} \Der(\log D)$), qui est localement libre de
rang 1. Pour des raison g\'en\'erales\footnote{On utilise le fait que
$\DX$ (resp. $\VCERO$) est l'alg\`ebre enveloppante de l'alg\'ebro\"{\i}de
de Lie $\Der_{\C_X}(\OX)$ (resp $\Der(\log D)$; voir
\cite{calde_ens}, prop. 2.2.5)} (cf.
\cite{rine-63,MR2000f:53109,MR2000m:32015}), $\omega_X$ (resp.
$\omega_X(\log D)$) a une structure naturelle de $\DX$-module \`a
droite (resp. de $\VCERO$-module \`a droite) et l'inclusion naturelle
$\omega_X \subset \omega_X(\log D)$ est $\VCERO$-lin\'eaire.
Rappelons que l'action \`a droite d'un champ de vecteurs $\delta$ sur
une $n$-forme diff\'erentielle $\theta$ est donn\'ee par $\theta
\delta = -L_{\delta}(\theta)$, o\`u $L_{\delta}$ est la d\'eriv\'ee
de Lie.

\section{Op\'erations externes}

Si $\M,\N$ sont deux $\VCERO$-modules \`a gauche, alors les
$\OX$-modules $\Hom_{\OX}(\M,\N)$ et $\M\otimes_{\OX} \N$ ont une
structure naturelle de $\VCERO$-module \`a gauche: $$ (\delta h)(m) =
- h(\delta m) + \delta h(m), \quad \delta \left( m\otimes n\right) =
(\delta m)\otimes n + m\otimes (\delta n)$$o\`u $\delta$ est une
d\'erivation logarithmique, $h$ est une section locale de
$\Hom_{\OX}(\M,\N)$ et $m, n$ sont des sections locales de $\M,\N$
respectivement.

On voit facilement que les isomorphismes naturels $\OX((m+m')D)
\simeq \OX(mD)\otimes_{\OX} \OX(m'D) $ et
$\Hom_{\OX}(\OX(mD),\OX)\simeq \OX(-mD)$ sont $\VCERO$-lin\'eaires.

Si $\E$ est une connexion logarithmique int\'egrable,
$\Hom_{\OX}(\E,\OX)$ est aussi une connexion logarithmique
int\'egrable qu'on notera $\E^*$. Dans ce cas, si $\M$ est un autre
$\VCERO$-module \`a gauche, l'isomorphisme canonique de $\OX$-modules
\begin{equation}\label{eq:can1}
\E^*\otimes_{\OX} \M \xrightarrow{} \Hom_{\OX}(\E,\M)
\end{equation}est un isomorphisme de $\VCERO$-modules \`a gauche.

Si $\PP, \QQ$ sont deux $\VCERO$-modules \`a droite, alors le
$\OX$-module $\Hom_{\OX}(\PP,\QQ)$ a une structure naturelle de
$\VCERO$-module \`a gauche: $$ (\delta h)(p) = h(p\delta) -
h(p)\delta,
$$o\`u $\delta$ est une d\'erivation logarithmique,
$h$ est une section locale de $\Hom_{\OX}(\PP,\QQ)$ et $p$ est une
section locale de $\PP$.

Si $\PP$ (resp. $\N$) est un $\VCERO$-module \`a droite (resp \`a
gauche), alors le $\OX$-module $\PP\otimes_{\OX} \N$ a une structure
naturelle de $\VCERO$-module \`a droite: $$ \left(p\otimes n\right)
\delta = (p\delta)\otimes n - p\otimes (\delta n)$$o\`u $\delta$ est
une d\'erivation logarithmique et $p, n$ sont des section locales de
$\PP,\N$ respectivement.

Le lemme suivant se d\'emontre comme les assertions (A.4) et (A.6)
dans \cite{es_vi_86}:

\begin{lema} \label{A.4} Si $\M,\N$ sont deux $\VCERO$-modules \`a gauche, on a un
isomorphisme naturel $\C_X$-lin\'eaire
$$ \Hom_{\VCERO}(\M,\N) \simeq \Hom_{\VCERO}(\OX,
\Hom_{\OX}(\M,\N)).$$Si de plus $\M$ est une connexion logarithmique
int\'egrable on a un isomorphisme naturelle dans la cat\'egorie
d\'eriv\'ee
$$ R\Hom_{\VCERO}(\M,\N) \simeq R\Hom_{\VCERO}(\OX,\M^*\otimes_{\OX}
\N).$$
\end{lema}

\numero {\it (Commutativit\'e et associativit\'e du produit tensoriel
externe)} \label{nume-1}
\smallskip

Si $\M,\N$ sont deux $\VCERO$-modules \`a gauche, l'isomorphisme
canonique de $\OX$-modules $\M \otimes_{\OX} \N\simeq \N
\otimes_{\OX} \M$ est en fait un isomorphisme de $\VCERO$-modules \`a
gauche.

Si $\PP$ est un autre $\VCERO$-module \`a gauche (resp. \`a droite)
l'isomorphisme canonique de $\OX$-modules
$$ \left( \PP\otimes_{\OX} \M\right) \otimes_{\OX} \N \simeq
\PP\otimes_{\OX} \left( \M \otimes_{\OX} \N\right)$$est
$\VCERO$-lin\'eaire \`a gauche (resp. \`a droite).
\medskip

D'apr\`es \cite{ksaito_log}, on a $\omega_X(\log D) = \omega_X
\otimes_{\OX} \OX(D)$.
\medskip

La preuve de la proposition suivante est facile.

\begin{nproposicion} \label{prop-ome} La structure naturelle de $\VCERO$-module \`a
droite sur $\omega_X(\log D)$ co\"{\i}ncide avec celle de $\omega_X
\otimes_{\OX} \OX(D)$ provenant de la structure naturelle de
$\DX$-module \`a droite sur $\omega_X$, et donc de $\VCERO$-module
\`a droite, et de la structure naturelle de $\VCERO$-module \`a
gauche sur $\OX(D)$.
\end{nproposicion}

Si $\M$ est un $\VCERO$-module \`a gauche, le $\OX$-module
$\M\otimes_{\OX} \VCERO$ est un $\VCERO$-bimodule: la structure \`a
gauche est celle qu'on vient de d\'efinir \`a partir de celles de de
$\M$ et de $\VCERO$, tandis que la structure \`a droite provient
seulement de celle de $\VCERO$. De fa\c{c}on analogue, le $\OX$-module
$\VCERO\otimes_{\OX} \M$ est aussi un $\VCERO$-bimodule: la structure
\`a droite provient de la structure \`a droite de $\VCERO$ et de la
structure \`a gauche de $\M$, tandis que la structure \`a gauche
provient seulement de celle de $\VCERO$.

\begin{nlema} \label{lema-1} Sous les hypoth\`eses pr\'ec\'edentes, il existe un
isomorphisme naturel unique de $\VCERO$-bimodules $$
\VCERO\otimes_{\OX} \M \simeq \M\otimes_{\OX}\VCERO$$qui applique
$1\otimes m$ dans $m\otimes 1$.
\end{nlema}

\begin{prueba} Le morphisme canonique $m\in \M \mapsto m\otimes 1\in
\M\otimes_{\OX}\VCERO$ est $\OX$-lin\'eaire \`a gauche et s'\'etend
donc \`a un morphisme $\VCERO$-lin\'eaire \`a gauche $$\varphi:
P\otimes m\in \VCERO\otimes_{\OX} \M \mapsto P\cdot (m\otimes 1)\in
\M\otimes_{\OX}\VCERO.$$On d\'emontre facilement que $\varphi$ est
aussi $\VCERO$-lin\'eaire \`a droite.

De fa\c{c}on analogue, le morphisme canonique
$$m\in \M \mapsto 1\otimes m\in
\VCERO\otimes_{\OX}\M$$est $\OX$-lin\'eaire \`a droite et s'\'etend
donc \`a un morphisme $\VCERO$-lin\'eaire \`a droite $$\psi: m\otimes
P\in \M\otimes_{\OX} \VCERO \mapsto (1\otimes m)\cdot P\in
\VCERO\otimes_{\OX}\M,$$qui devient aussi $\VCERO$-lin\'eaire \`a
gauche.

On a $(\psi \pcirc \varphi)(P\otimes m) = \psi(P\cdot (m\otimes 1)) =
P\cdot \psi (m\otimes 1) = P\cdot (1\otimes m) = P\otimes m$, d'o\`u
$\psi \pcirc \varphi$ est l'identit\'e. De la m\^eme fa\c{c}on
$\varphi\pcirc \psi$ est aussi l'identit\'e, ce qui prouve le lemme.
\end{prueba}

\begin{ncorolario} \label{coro-lema-1} Si $\M,\N$ sont deux $\VCERO$-modules \`a
gauche, on a un morphisme naturel de $\VCERO$-modules \`a droite $$
\Hom_{\VCERO}(\M,\VCERO)\otimes_{\OX} \N \to
\Hom_{\VCERO}(\M,\N\otimes_{\OX}\VCERO),$$ qui est un isomorphisme si
$\N$ est une connexion logarithmique int\'egrable.
\end{ncorolario}

\begin{nproposicion} \label{A.8} Si $\PP$ est un $\VCERO$-module \`a
droite et $\M,\N$ sont deux $\VCERO$-modules \`a gauche, le
morphisme
$$ p\otimes (m\otimes n) \in \PP\otimes_{\VCERO} \left( \M\otimes_{\OX} \N\right)
\mapsto (p\otimes m)\otimes n \in \left(\PP\otimes_{\OX}
\M\right)\otimes_{\VCERO} \N$$est bien d\'efini et est un
isomorphisme $\C_X$-lin\'eaire.
\end{nproposicion}

\begin{prueba} Elle est laiss\'ee au lecteur.
\end{prueba}

\numero \label{nume-bimo} Si $\PP$ est un $\VCERO$-module \`a droite,
le $\OX$-module $\PP\otimes_{\OX} \VCERO$ est un $\VCERO$-bimodule
droite-droite: la premi\`ere structure \`a droite est celle qui
provient de la  structure \`a droite de $\PP$ et de la structure \`a
gauche de $\VCERO$, tandis que la deuxi\`eme structure \`a droite
provient de la structure \`a droite de $\VCERO$.

La preuve du lemme suivant est analogue \`a celle du lemme
\ref{lema-1}: \smallskip

\begin{nlema} \label{lema-2} Sous les hypoth\`eses pr\'ec\'edentes, il existe
une involution $\OX$-lin\'eaire du module $\PP\otimes_{\OX} \VCERO$
qui interchange les deux structures de $\VCERO$-module \`a droite.
\end{nlema}

\begin{ncorolario} \label{coro-1} Sous les hypoth\`eses pr\'ec\'edentes,
si $\PP$ est localement libre de rang fini comme $\OX$-module,
 alors le module
$\PP\otimes_{\OX} \VCERO$ est localement libre de rang fini en tant
que $\VCERO$-module \`a droite pour la premi\`ere structure du
paragraphe \ref{nume-bimo}.
\end{ncorolario}

\begin{nteorema} \label{teor-fund} Soit $\PP$ un $\DX$-module \`a droite et $\N$ un
$\VCERO$-module \`a gauche. Alors, le morphisme naturel
$$ p\otimes n \in \PP\otimes_{\OX} \N\mapsto p\otimes(1\otimes n)\in
\PP\otimes_{\OX} \left(\DX\otimes_{\VCERO} \N\right)$$ est
$\VCERO$-lin\'eaire \`a droite et le morphisme induit
$$ \lambda: \left( \PP\otimes_{\OX} \N\right) \otimes_{\VCERO} \DX
\xrightarrow{} \PP\otimes_{\OX} \left(\DX\otimes_{\VCERO} \N\right)$$
est un isomorphisme de $\DX$-modules \`a droite.
\end{nteorema}

\begin{prueba} Pour d\'emontrer le th\'eor\`eme il suffit de proc\'eder fibre \`a fibre, et
ceci est une cons\'equence du corollaire \ref{A.2} pour $U_0 =
\D_{X,p}(\log D)$, $U=\D_{X,p}$ et $p\in X$.
\end{prueba}

\begin{ncorolario} \label{cor-fund} Soit $\PP$ un $\DX$-module \`a droite qui
est localement libre de rang fini sur $\OX$ et $\N$ un
$\VCERO$-module \`a gauche qui admet des r\'esolutions localement
libres\footnote{Ceci sera le cas si par exemple $\N$ est un
$\VCERO$-module coh\'erent.}. Alors on a un isomorphisme naturel dans
la cat\'egorie d\'eriv\'ee des $\DX$-modules \`a droite
$$ \left( \PP\otimes_{\OX} \N\right) \Lotimes_{\VCERO} \DX
\xrightarrow{} \PP\otimes_{\OX} \left(\DX\Lotimes_{\VCERO}
\N\right).$$
\end{ncorolario}

\begin{prueba} D'apr\`es les hypoth\`eses, on peut se r\'eduire au cas o\`u
$\N=\VCERO$ (comme $\VCERO$-module \`a gauche!). Dans ce cas, le
corollaire \ref{coro-1} nous dit que le $\VCERO$-module \`a droite
$\PP\otimes_{\OX} \VCERO$ est localement libre de rang fini, donc
acyclique pour le produit tensoriel sur $\VCERO$. On conclut en
appliquant le th\'eor\`eme \ref{teor-fund}.
\end{prueba}

Les r\'esultats de l'Appendice nous permettent aussi de montrer la
proposition suivante.

\begin{nproposicion} Pour toute connexion logarithmique int\'egrable
$\E$, il existe un isomorphisme canonique de complexes de
$\VCERO$-modules \`a gauche:
$$ \SPlog (\E) \simeq \SPlog \otimes_{\OX} \E$$
qui est compatible avec les augmentations vers $\E$.
\end{nproposicion}

\begin{prueba}  D'apr\`es le corollaire \ref{A.7}, pour chaque
degr\'e $k=0,\dots, n$ on a des isomorphismes canoniques de
$\VCERO$-modules \`a gauche: \begin{eqnarray*} & Sp^{k}_{\VCERO}
(\E) = \VCERO \otimes_{\OX} \left[
\stackrel{k}{\bigwedge}\derlogD\otimes_{\OX} \E \right] \simeq&\\
&\left[ \VCERO \otimes_{\OX} \stackrel{k}{\bigwedge}\derlogD\right]
\otimes_{\OX} \E = Sp^{k}_{\VCERO}  \otimes_{\OX}
\E.&\end{eqnarray*}On v\'erifie facilement que ces isomorphismes
commutent avec les diff\'erentielles $\varepsilon^\bullet$ et
$d^\bullet \otimes \Id_{\E}$ et les augmentations vers $\E$
d\'efinies dans la section \ref{sec:1}.
\end{prueba}

\begin{nnota}
Tous les r\'esultats de cette section sont valables pour un
alg\'ebro\"{\i}de de Lie arbitraire et son alg\`ebre enveloppante \`a la
place de $\Der(\log D)$ et $\VCERO$ respectivement.
\end{nnota}

\section{Dualit\'e sur les connexions logarithmiques et sur les complexe de de Rham logarithmiques}

La proposition suivante est bien connue dans le cas des
$\DX$-modules. Elle se g\'en\'eralise sans peine au cas des
alg\'ebro\"{\i}des de Lie (cf. \cite{MR2000m:32015}, prop. 3.2.1).
\medskip

\begin{proposicion} \label{prop-0} Pour toute connexion logarithmique int\'egrable $\E$ le complexe
$R\Hom_{\VCERO}(\E,\VCERO)$ est concentr\'e en degr\'e $n$ et on a un
isomorphisme naturel de $\VCERO$-modules \`a droite
$$\Ext^n_{\VCERO}(\E,\VCERO) \simeq \omega_X(\log
D)\otimes_{\OX} \E^*.$$
\end{proposicion}

\begin{prueba} D'apr\`es le lemme \ref{A.4} et le corollaire
\ref{coro-lema-1} on a des isomorphismes naturels
$\VCERO$-lin\'eaires \`a droite: \begin{eqnarray*} &
R\Hom_{\VCERO}(\E,\VCERO) \simeq
R\Hom_{\VCERO}(\OX, \E^*\otimes_{\OX}\VCERO)&\\
&\simeq R\Hom_{\VCERO}(\OX, \VCERO)\otimes_{\OX}\E^*.&
\end{eqnarray*}On est donc r\'eduit au cas $\E=\OX$. Consid\'erons la
r\'esolution $\SP_{\VCERO}$ de $\OX$: \begin{eqnarray*} &
R\Hom_{\VCERO}(\OX, \VCERO)&\\ &= \Hom_{\VCERO}(\SP_{\VCERO},\VCERO)=
\Omega^{\bullet}_X(\log D)(\VCERO)&\end{eqnarray*}et l'augmentation
$$\Omega^n_X(\log D)\otimes_{\OX} \VCERO = \omega_X(\log D) \otimes_{\OX}
\VCERO\to \omega_X(\log D)$$donn\'ee par $\theta\otimes P \mapsto
\theta\cdot P$. Filtrons le complexe augment\'e par
$$ F^i \left(\Omega^k_X(\log D)\otimes_{\OX} \VCERO\right) =
\Omega^k_X(\log D)\otimes_{\OX} F^{i+k}\VCERO$$ et $F^i\omega_X(\log
D) = \omega_X(\log D)$ pour $k=0,\dots, n$ et $i\in\ZZ$. On v\'erifie
sans peine que le gradu\'e associ\'e est exact (cf. \cite{bou_a_10},
\S 9, 3), d'o\`u le r\'esultat.
\end{prueba}
\medskip

\setcounter{numero}{0} \numero \label{dual-V} Rappelons que le
foncteur de dualit\'e au niveau de $\DX$-modules, $\Dual_{\DX}:
\cd{b}{coh}{\DX}\to \cd{b}{coh}{\DX}$, est d\'efini par:
$$\Dual_{\DX} \M = \Hom_{\OX}(\omega_X,R\Hom_{\DX}(\M,\DX))[n].$$

De fa\c{c}on analogue nous pouvons consid\'erer le foncteur (cf.
\cite{MR2000m:32015}, 3.2)
$$\Dual_{\VCERO}: \cd{b}{coh}{\VCERO} \to \cd{b}{coh}{\VCERO}$$d\'efini
par
$$\Dual_{\VCERO}\M = \Hom_{\OX}(\omega_X(\log
D),R\Hom_{\VCERO}(\M,\VCERO))[n].$$ Il s'agit d'une involution de la
cat\'egorie triangul\'ee $\cd{b}{coh}{\VCERO}$.

La proposition pr\'ec\'edente nous dit que si $\E$ est une connexion
logarithmique int\'egrable  alors $\Dual_{\VCERO}(\E) = \E^*$.
\medskip

Si $\E$ est une connexion logarithmique int\'egrable, on notera
$\E(mD)$ la connexion logarithmique int\'egrable  $\E\otimes_{\OX}
\OX(mD)$.

\begin{nteorema} \label{main} Pour toute connexion logarithmique
int\'egrable $\E$ on a un isomorphisme naturel dans la cat\'egorie
d\'eriv\'ee des $\DX$-modules \`a droite: $$ R\Hom_{\DX}(
\DX\Lotimes_{\VCERO} \E,\DX) \simeq \omega_X\otimes_{\OX}
\left(\DX\Lotimes_{\VCERO} \E^*(D)\right)[-n].$$
\end{nteorema}

\begin{prueba} On a un isomorphisme canonique
$$ R\Hom_{\DX}( \DX\Lotimes_{\VCERO} \E,\DX) \simeq R\Hom_{\VCERO}(
\E,\DX),$$et par coh\'erence de $\VCERO$ et $\E$
$$ R\Hom_{\VCERO}(\E,\DX)\simeq
R\Hom_{\VCERO}(\E,\VCERO)\Lotimes_{\VCERO} \DX.$$D'apr\`es la
proposition \ref{prop-0} on a
$$ R\Hom_{\DX}( \DX\Lotimes_{\VCERO} \E,\DX) \simeq \left[\omega_X(\log
D)\otimes_{\OX} \E^*\right]\Lotimes_{\VCERO} \DX [-n].$$Pour conclure
il suffit d'appliquer \ref{nume-1}, la proposition \ref{prop-ome} et
le corollaire \ref{cor-fund}.
\end{prueba}

\begin{ncorolario} \label{cor-main} Pour toute connexion logarithmique
int\'egrable $\E$ on a un isomorphisme naturel dans la cat\'egorie
d\'eriv\'ee des $\DX$-modules \`a gauche:
$$ \Dual_{\DX}(\DX\Lotimes_{\VCERO} \E)\simeq \DX\Lotimes_{\VCERO}
\E^*(D).$$
\end{ncorolario}

\begin{nnota} \label{}  Dans le cas d'un diviseur de
Spencer et $\E =\OX$, le corollaire pr\'ec\'edent est une version
intrins\`eque du th\'eor\`eme de dualit\'e de \cite{cas_ucha_stek},
th. 4.3 (voir la pr\'esentation (\ref{eq:presenta})).
\end{nnota}

\begin{ncorolario}\label{dual-admi} Si $\E$ est une connexion logarithmique int\'egrable
admissible (voir d\'ef. \ref{def-admi}), alors $\E^*(D)$ est aussi
admissible. \end{ncorolario}

\begin{prueba} Comme le complexe $\DX\Lotimes_{\VCERO} \E $ est
\`a cohomologie holonome concentr\'ee en degr\'e $0$, son dual au
sens de $\DX$-modules l'est aussi et donc, d'apr\`es le corollaire
pr\'ec\'edent, $\E^*(D)$ est admissible.\end{prueba}

\begin{ncorolario} \label{cor:DRlog-DRusual} Pour toute connexion logarithmique
int\'egrable $\E$ on a un isomorphisme naturel dans la cat\'egorie
d\'eriv\'ee
$$ \Omega^{\bullet}_X(\log D)(\E)  \simeq R \Hom_{\DX}(\OX, \DX\Lotimes_{\VCERO}
\E(D)).$$
\end{ncorolario}

\begin{prueba} D'apr\`es (\ref{eq:calde-ens-3.2.1}), $\Omega_X^{\bullet}(\log D)(\E)
\simeq R\Hom_{\VCERO}(\OX,\E)$. Or,
\begin{eqnarray*} & R\Hom_{\VCERO}(\OX,\E) \simeq
R\Hom_{\VCERO}(\Dual_{\VCERO}\E, \Dual_{\VCERO}\OX)\simeq &\\
&\simeq R\Hom_{\VCERO}(\E^*, \OX) \simeq
R\Hom_{\DX}(\DX\Lotimes_{\VCERO} \E^*, \OX)\simeq
&\\
& R\Hom_{\DX}(\Dual_{\DX}\OX, \Dual_{\DX}(\DX\Lotimes_{\VCERO}
\E^*)) \simeq&\\ & R\Hom_{\DX}(\OX, \DX\Lotimes_{\VCERO} \E(D)).&
\end{eqnarray*}
\end{prueba}

\begin{ncorolario} \label{cor:admi-equiv} Soit $\E$ une connexion logarithmique int\'egrable.
Les propri\'et\'es suivantes sont \'equivalentes:
\begin{enumerate}
\item[1)] $\E^*$ est admissible.
\item[2)] $\E(D)$ est admissible.
\item[3)] Le complexe de de Rham logarithmique $\Omega^{\bullet}_X(\log D)(\E)$
est un faisceau pervers.
\end{enumerate}
\end{ncorolario}

\begin{prueba} L'\'equivalence entre 1) et 2) a \'et\'e prouv\'ee
dans le corollaire \ref{dual-admi}.
\smallskip

\noindent 2) $\Rightarrow$ 3) est une cons\'equence directe de la
d\'efinition des connexions logarithmiques admissibles, du corollaire
\ref{cor:DRlog-DRusual} et du fait que le complexe de de Rham (usuel)
d'un $\DX$-module holonome est un faisceau pervers.
\smallskip

\noindent 3) $\Rightarrow$ 1):\ Au cours de la preuve du corollaire
\ref{cor:DRlog-DRusual} on a exprim\'e le complexe $\K =
\Omega^{\bullet}_X(\log D)(\E)$ comme le complexe des solutions
holomorphes du complexe born\'e de $\DX$-modules \`a cohomologie
coh\'erente $\M = \DX\Lotimes_{\VCERO} \E^*$. D'apr\`es
\cite{meb_formalisme}, chap. II, th. (4.1.5), et la
constructibilit\'e de $\K$ on d\'eduit d'abord que le complexe $\M$
est \`a cohomologie holonome. Or, par le th\'eor\`eme de bidualit\'e
locale \cite{meb_formalisme}, chap. I, th. (10.13), on sait que
$$\DX^{\infty}\otimes_{\DX} \M \simeq R \Hom_{\C_X}(\K, \OX),$$et
comme $\K$ est pervers, le deuxi\`eme complexe est, toujours par le
th\'eor\`eme de bidualit\'e locale, concentr\'e en degr\'e $0$. Pour
conclure il suffit d'invoquer la fid\`ele platitude de $\DX^{\infty}$
sur $\DX$ \cite{skk} (voir aussi \cite{nar-rojas}).
\end{prueba}

\begin{ncorolario} \label{cor:crit-spencer} Les propri\'et\'es suivants sont \'equivalentes:
\begin{enumerate}
\item Le diviseur $D$ est de Spencer (i.e. $\OX$ est
admissible).
\item Le complexe de de Rham logarithmique $\Omega^{\bullet}_X(\log
D)$ est un faisceau pervers.
\end{enumerate}
\end{ncorolario}

Le corollaire suivant g\'en\'eralise la proposition (A.2) de
\cite{es_vi_86}. Elle s'applique \`a toutes les connexions
logarithmiques int\'egrables par rapport aux croisements normaux, et
plus g\'en\'eralement par rapport aux diviseurs libres localement
quasi-homog\`enes \cite{calde_nar_LQHKF}.

\begin{ncorolario} \label{cor:dual-DRlog} Soit $\E$ une connexion logarithmique int\'egrable  (par rapport \`a
$D$) telle que $\E^*$ est admissible (et donc $\E(D)$ est aussi
admissible). Alors on a un isomorphisme naturel dans la cat\'egorie
d\'eriv\'ee
$$ \Omega_X^{\bullet}(\log D)(\E) \simeq \Omega_X^{\bullet}(\log
D)(\E^*(-D))^{\vee},$$o\`u $\vee$ d\'enote le dual de Verdier.
\end{ncorolario}

\begin{prueba} D'apr\`es le corollaire \ref{cor:DRlog-DRusual},
on a des isomorphismes naturels
\begin{eqnarray*}
& \Omega^{\bullet}_X(\log D)(\E)  \simeq R \Hom_{\DX}(\OX,
\DX\Lotimes_{\VCERO} \E(D)),&\\
&\Omega_X^{\bullet}(\log D)(\E^*(-D)) \simeq R \Hom_{\DX}(\OX,
\DX\Lotimes_{\VCERO} \E^*).& \end{eqnarray*} Le corollaire est donc
une cons\'equence du th\'eor\`eme de dualit\'e locale (cf.
\cite{meb_formalisme}, ch. I, th. (4.3.1); voir aussi \cite{nar-ldt})
et du corollaire \ref{cor-main}.
\end{prueba}

\begin{nnota} D'apr\`es la proposition \ref{K-admi}, le corollaire pr\'ec\'edent
s'applique \`a toute connexion logarithmique int\'egrable dans le
cas o\`u le diviseur $D$  est de Koszul.
\end{nnota}

\section{Un crit\`ere diff\'erentiel pour le th\'eor\`eme de comparaison logarithmique}

Notons $\rho: \DX\Lotimes_{\VCERO} \OX(D)\to \OX(\star D)$ le
morphisme $\DX$-lin\'eaire \`a gauche donn\'e par $\rho(P\otimes a) =
P(a)$.

Rappelons qu'on dit que {\it le th\'eor\`eme de comparaison
logarithmique} est vrai pour $D$ si l'inclusion
$$ \Omega_X^{\bullet}(\log D) \hookrightarrow \Omega_X^{\bullet}(\star
D)$$est un quasi-isomorphisme, ce qui \'equivaut, gr\^ace au
th\'eor\`eme de comparaison de Grothendieck, au fait que le morphisme
naturel
$$\Omega_X^{\bullet}(\log D) \to R j_* j^{-1} \Omega^{\bullet}_X$$est un quasi-isomorphisme
(voir \cite{cas_mond_nar_96,calde_mond_nar_cas}), o\`u $j: X-D
\hookrightarrow X$ est l'inclusion ouverte.
\medskip

Plus g\'en\'eralement, pour chaque connexion logarithmique
int\'egrable $\E$ on consid\`ere le morphisme
$$\rho_{\E}: \DX\Lotimes_{\VCERO} \E(D)\to \E(\star D),$$ donn\' e
par $\rho_{\E}(P\otimes e') = P(e')$, o\`u $\E(\star D)$ est le
$\DX$-module obtenu par restriction de scalaires \`a partir du
$\VCERO(\star D) =\DX(\star D)$-module $\OX(\star D)\otimes_{\OX}
\E$. Nous savons que $\E(\star D)$ est une connexion m\'eromorphe et
donc holonome (cf. \cite{meb_nar_dmw} th. 4.1.3). En fait, $\E(\star
D)$ est r\'eguli\`ere dans la partie lisse de $D$ (elle a des p\^oles
logarithmiques!), et donc elle est r\'eguli\`ere partout
\cite{meb-cimpa-2}, cor. 4.3-14.
\medskip

Aussi, le complexe de de Rham logarithmique $\Omega_X^{\bullet}(\log
D)(\E)$ est un sous-complexe du complexe de de Rham m\'eromorphe
$\Omega_X^{\bullet}(\log D)(\E(\star D))=\Omega_X^{\bullet}(\E(\star
D))$. Les restrictions \`a $X-D$ des ces deux complexes co\"{\i}ncident et
sont quasi-isomorphes au syst\`eme local $\L$ des section
horizontales de $\E$ sur $X-D$.

\begin{teorema} \label{teor:gordo} Sous les hypoth\`eses pr\'ec\'edentes,
les propri\'et\'es suivantes sont \'equivalentes:
\begin{enumerate}
\item[1)] Le morphisme canonique $\Omega_X^{\bullet}(\log D)(\E)
\to R j_* \L$ est un isomorphisme dans la cat\'egorie d\'eriv\'ee.
\item[2)] L'inclusion $\Omega_X^{\bullet}(\log D)(\E) \hookrightarrow \Omega_X^{\bullet}(\E(\star
D))$ est un quasi-isomorphisme.
\item[3)] Le morphisme $\rho_{\E}: \DX\Lotimes_{\VCERO} \E(D)\to \E(\star D)$
est un isomorphisme dans la cat\'egorie d\'eriv\'ee.
\end{enumerate}
\end{teorema}

\begin{prueba} Notons pour simplifier $\VO = \DX(\log D)$.
L'\'equivalence entre 1) et 2) provient de la r\'egularit\'e de
$\E(\star D)$:
$$ \Omega_X^{\bullet}(\E(\star
D)) \simeq R j_* j^{-1} \Omega_X^{\bullet}(\E(\star D)) \simeq R
j_* \L.$$
\smallskip

\noindent 3) $\Rightarrow$ 2):\ Tout consiste \`a consid\'erer le
diagramme commutatif suivant dans la cat\'egorie d\'eriv\'ee des
complexes de faisceaux de $\C$-espaces vectoriels:
\begin{equation*}
\begin{CD}
\Omega_X^{\bullet}(\log D)(\E) @>{\text{incl.}}>>
\Omega_X^{\bullet}(\log
D)(\E(\star D))\\
@V{\simeq}VV    @V{\simeq}VV\\
\Omega_X^{\bullet}(\log D)(\VO)\otimes_{\VO} \E  @>{\Id \otimes
\text{incl.}}>> \Omega_X^{\bullet}(\log D)(\VO)\otimes_{\VO}
\E(\star D)\\
@V{\simeq}VV    @V{\simeq}VV\\
\left(\omega_X(\log D) \Lotimes_{\VO} \E\right)[-n] @>{\Id \otimes
\text{incl.}}>> \left(\omega_X(\log D) \Lotimes_{\VO} \E(\star
D)\right)[-n]\\
@V{\simeq}V{(*)}V    @V{\simeq}V{(*)}V\\
\left(\omega_X \Lotimes_{\VO} \E(D)\right)[-n] @>{\Id \otimes
\text{incl.}}>> \left(\omega_X \Lotimes_{\VO} \E(\star
D)\right)[-n]\\
@V{\simeq}VV    @V{\simeq}VV\\
\left(\omega_X \Lotimes_{\DX} \DX \Lotimes_{\VO} \E(D)\right)[-n]
@>{\Id\otimes \Id \otimes \text{incl.} }>> \left(\omega_X
\Lotimes_{\DX} \DX \Lotimes_{\VO} \E(\star
D)\right)[-n]\\
@V{=}VV    @V{\simeq}V{(**)}V\\
\omega_X[-n]\Lotimes_{\DX} \left(\DX \Lotimes_{\VO} \E(D)\right)
@>{\Id\otimes \rho_{\E}}>> \omega_X[-n]\Lotimes_{\DX} \E(\star
D)\\
@V{\simeq}VV    @V{\simeq}VV\\
\DR (\DX \Lotimes_{\VO} \E(D)) @>{\DR (\rho_{\E})}>> \DR (\E(\star
D)),
\end{CD}
\end{equation*}
o\`u les isomorphismes $(*)$ sont ceux qui proviennent de la
proposition \ref{A.8}, et l'isomorphisme $(**)$ est donn\'e par
$$ \textstyle \DX \Lotimes_{\VO} \E (\star D) \simeq \DX \Lotimes_{\VO}
\VO(\star D) \otimes_{\VO(\star D)} \E(\star D) \simeq \VO(\star
D) \otimes_{\VO(\star D)} \E(\star D) \simeq \E(\star D).$$ Comme
les fl\`eches verticales sont des isomorphismes, si $\rho_{\E}$
est un isomorphisme dans la cat\'egorie d\'eriv\'ee, alors
l'inclusion $\Omega_X^{\bullet}(\log D)(\E) \hookrightarrow
\Omega_X^{\bullet}(\log D)(\E(\star D))$ est un
quasi-isomorphisme.
\smallskip

\noindent 2) $\Rightarrow$ 3):\ R\'eciproquement, par le diagramme
pr\'ec\'edent, la propri\'et\'e 2) entra\^{\i}ne que $\DR
(\rho_{\E})$ est un isomorphisme dans la cat\'egorie d\'eriv\'ee.
Soit $\QQ$ le c\^one de $\rho_{\E}$. Il s'agit d'un complexe de
$\DX$-modules \`a cohomologie coh\'erente dont le de Rham est nul.
Par l'argument de la preuve de 3) $\Rightarrow$ 1) dans le corollaire
\ref{cor:admi-equiv} on d\'eduit que $\QQ=0$ et donc $\rho_{\E}$ est
un isomorphisme dans la cat\'egorie d\'eriv\'ee.
\end{prueba}

\begin{corolario} \label{criterio} Les propri\'et\'es suivantes sont \'equivalentes:
\begin{enumerate}
\item Le th\'eor\`eme de comparaison logarithmique est vrai pour $D$.
\item Le morphisme $\rho: \DX\Lotimes_{\VCERO} \OX(D)\to \OX(\star D)$
est un isomorphisme dans la cat\'egorie d\'eriv\'ee.
\end{enumerate}
\end{corolario}

Le corollaire \ref{criterio} nous dit en particulier qu'une condition
n\'ecessaire pour que le th\'eor\`eme de comparaison logarithmique
soit vrai pour $D$ est que la connexion logarithmique int\'egrable
$\OX(D)$ soit admissible, ou encore par le corollaire
\ref{dual-admi}, que $\OX$ soit admissible, ce qui revient \`a dire
que $D$ soit un diviseur de Spencer dans la terminologie de
\cite{cas_ucha_stek}, def. 3.3.
\medskip

Si $f=0$ est une \'equation r\'eduite locale de $D$,
$\delta_1,\dots,\delta_n$ est une base locale de $\derlogD$ et
$\delta_i(f)=\alpha_i f$, en tenant compte de la pr\'esentation
(\ref{eq:presenta}), le th\'eor\`eme \ref{criterio} peut s'\'enoncer
comme l'\'equivalence des propri\'et\'es suivantes:
\begin{enumerate}
\item[a)] Le th\'eor\`eme de comparaison logarithmique est vrai pour $D$.
\item[b)] (b-1) $D$ est un diviseur de Spencer, (b-2) le $\DX$-module des fonctions
m\'eromorphes le long de $D$ est engendr\'e par $f^{-1}$ et (b-3)
l'annulateur de $f^{-1}$ est le $\DX$-id\'eal \`a gauche engendr\'e
par $\delta_1+\alpha_1,\dots, \delta_n+\alpha_n$.
\end{enumerate}
En fait, gr\^ace \`a la proposition 1.3 de \cite{torre-45-bis}, la
condition (b-3) entra\^{\i}ne la condition (b-2). Par cons\'equent,
on a le corollaire suivant:

\begin{corolario} \label{coro-torre}  Les propri\'et\'es suivantes sont \'equivalentes:
\begin{enumerate}
\item Le th\'eor\`eme de comparaison logarithmique est vrai pour $D$.
\item Le complexe $\DX\Lotimes_{\VCERO} \OX(D)$ est concentr\'e en
degr\'e $0$ et le morphisme $\rho: \DX\otimes_{\VCERO} \OX(D)\to
\OX(\star D)$ est injectif.
\end{enumerate}
\end{corolario}

En utilisant le corollaire \ref{criterio} et l'argument de
\cite{cas_ucha_stek}, th. 5.2 et lemme 5.3, nous allons donner une
nouvelle preuve diff\'erentielle non topologique du th\'eor\`eme de
comparaison logarithmique pour les diviseurs libres localement
quasi-homog\`enes.

\begin{teorema} \label{tcl}(\cite{cas_mond_nar_96}) Si $D$ est un diviseur libre
localement quasi-homog\`ene, alors l'inclusion
$$\Omega_X^{\bullet}(\log D) \hookrightarrow \Omega_X^{\bullet}(\star
D)$$est un quasi-isomorphisme.
\end{teorema}

\begin{prueba} D'apr\`es \cite{calde_nar_LQHKF} et la proposition \ref{K-admi}
nous savons que le complexe $\DX\Lotimes_{\VCERO}\OX(mD)$ est
holonome et concentr\'e en degr\'e $0$ pour tout $m\geq 1$.

Tout d'abord nous allons calculer le cycle caract\'eristique des
$\DX\otimes_{\VCERO}\OX(mD)$. Soit $f=0$ une \'equation locale
r\'eduite de $D$ au voisinage d'un point $p\in D$ et
$\delta_1,\dots,\delta_n$ une base de $\Der(\log D)_p$, avec
$\delta_i(f)=\alpha_i f$. Le module $\OO_{X,p}(mD)$, $m\geq 1$,
 est engendr\'e comme $\D_{X,p}(\log D)$-module par $f^{-m}$, dont l'annulateur est
 l'id\'eal \`a gauche $I_m$ engendr\'e par
 $\delta_1+m\alpha_1,\dots,\delta_n+m\alpha_n$ (voir
 (\ref{eq:presenta})).

 D'apr\`es \cite{calde_nar_compo}, Remark 5.10,
 $\D_{X,p}/\D_{X,p} I_m \simeq \OO_{X,p}(\star D)$ pour $m\gg 0$ (il faut que $-m$ soit plus
 petit que la plus petite racine enti\`ere
 du polyn\^ome de Bernstein-Sato de $f$). Or,
 comme les symboles $\sigma(\delta_i+m\alpha_i)=\sigma(\delta_i)$,
 $i=1,\dots, n$ forment une suite r\'eguli\`ere dans le gradu\'e de $\D_{X,p}$
 pour tout $m\geq 0$ (le diviseur $D$ est de Koszul, d'apr\`es \cite{calde_nar_LQHKF}), nous
 d\'eduisons que le cycle caract\'eristique de tous les $\D_{X,p}/\D_{X,p}
 I_m$ co\"{\i}ncide avec le cycle caract\'eristique de $\OO_{X,p}(\star
 D)$.

 Les raisonnements locaux pr\'ec\'edents se recollent et on conclut que
 tous les $\DX$-modules holonomes $\DX\otimes_{\VCERO}\OX(mD)$,
 $m\geq 1$, ont m\^eme cycle caract\'eristique \'egal au cycle
 caract\'eristique de $\OX(\star D)$.

D'apr\`es le corollaire \ref{criterio}, pour d\'emontrer le
th\'eor\`eme de comparaison logarithmique pour $D$ il suffit de
d\'emontrer que le morphisme $\DX$-lin\'eaire:
\begin{equation} \label{eq-induc}
\rho_{D,X}: P\otimes a \in \DX\otimes_{\VCERO}\OX(D)\mapsto
P(a)\in\OX(\star D)
\end{equation}
est un isomorphisme, ou m\^eme qu'il est injectif car les deux
modules ont m\^eme cycle caract\'eristique.

Nous allons proc\'eder par r\'ecurrence sur la dimension $n$ de la
vari\'et\'e ambiante $X$ en utilisant la m\'ethode\footnote{Cette
m\'ethode a \'et\'e utilis\'ee aussi dans \cite{calde_nar_LQHKF}, th.
3.2; \cite{cas_ucha_stek}, th. 5.2; \cite{calde_nar_compo}, prop.
5.2, th. 5.9.} de \cite{cas_mond_nar_96}, prop. 2.4.

Si $n=\dim X = 1$, le r\'esultat est clair. Supposons  que le
morphisme (\ref{eq-induc}) est un isomorphisme si $D$ est un diviseur
libre localement quasi-homog\`ene dans une vari\'et\'e analytique
complexe lisse $X$ de dimension $n-1$.

Soit maintenant $D\subset X$ un diviseur libre localement
quasi-homog\`ene dans une vari\'et\'e analytique complexe lisse $X$
de dimension $n$ et $p\in D$. D'apr\`es \cite{cas_mond_nar_96}, prop.
2.4 et lemme 2.2,~(iv), il existe un voisinage ouvert $U$ de $p$ tel
que pour tout $q\in U\cap D$, $q\neq p$, le germe $(X,D,q)$ est
analytiquement isomorphe  \`a un produit $ (\C^{n-1}\times \C,
D^{\prime}\times \C, (0,0))$, avec $D^{\prime}\subset \C^{n-1}$
diviseur libre et localement quasi-homog\`ene.

Posons $X'=\C^{n-1}$ et notons $\pi:X'\times \C\to X'$ la projection.
Par l'hypoth\`ese de r\'ecurrence, le morphisme
$$\rho_{D',X'}:\D_{X'}\otimes_{\D_{X'}(\log D')} \OO_{X'}(D') \to
\OO_{X'}(\star D')$$est un isomorphisme. Or, le morphisme
$$\rho_{D'\times\C,X'\times\C}:
\D_{X'\times\C}\otimes_{\D_{X'\times\C}(\log D'\times\C)}
\OO_{X'\times\C}(D'\times\C)\to \OO_{X'\times\C}(\star (
D'\times\C))$$s'identifie avec $\pi^*\rho_{D',X'}$ et donc est un
isomorphisme. Par cons\'equent, le morphisme
$$\rho_{D,X}: \DX\otimes_{\VCERO}\OX(D)\xrightarrow{}\OX(\star
D)$$ est un isomorphisme sur $U-\{p\}$.

Notons $\K$ le noyau de la restriction \`a $U$ de $\rho_{D,X}$. Il
s'agit d'un module holonome support\'e par $p$. Par l'argument de
\cite{cas_ucha_stek}, th. 5.2 et lemme 5.3, nous d\'eduisons que $\K$
est nul, et par cons\'equent $\rho_{D,X}$ est injective, ce qui
termine la preuve du th\'eor\`eme.
\end{prueba}
\medskip

\begin{corolario} Si $D$ est un diviseur libre localement
quasi-homog\`ene et $j:X-D\hookrightarrow X$ est l'inclusion, on a
des isomorphismes canoniques:
$$\Omega_X^{\bullet}(\log D) \simeq R j_* \C_{X-D},\quad j_!
\C_{X-D} \simeq \Omega_X^{\bullet}(\log D)(\OX(-D)).$$
\end{corolario}

\begin{prueba} Le premier isomorphisme est une cons\'equence du
th\'eor\`eme \ref{tcl} et du th\'eor\`eme de comparaison de
Grothendieck $\Omega_X^{\bullet}(\star D) \simeq R j_* \C_{X-D}$.

Le deuxi\`eme isomorphisme est une cons\'equence du premier et du
corollaire \ref{cor:dual-DRlog}.
\end{prueba}

\begin{nota} Notons que dans le corollaire pr\'ec\'edent, le
morphisme canonique $j_! \C_{X-D} \to R j_* \C_{X-D}$ correspond \`a
l'inclusion $\OX(-D) \hookrightarrow \OX$.
\end{nota}

\section{Exemples et questions}

L'exemple suivant r\'epond au probl\`eme 6.6. de
\cite{calde_nar_compo} (voir \ref{cor:crit-spencer}).

\begin{ejemplo}\footnote{Dans cet exemple nous avons utilis\'e \cite{M2} et \cite{D-mod-M2}.} Soit $D\subset X=\C^3$ le diviseur donn\'e par
l'\'equation r\'eduite
$$h=(xz+y)(x^4+y^5+xy^4)=x^2y^4z+xy^5z+xy^5+y^6+x^5z+x^4y=0, $$qui
a \'et\'e consid\'er\'e dans \cite{cas_ucha_pams}. Le diviseur $D$
est libre, une base de $\Der(\log D)$ \'etant $\{\delta_1,
\delta_2, \delta_3 \}$, avec  {\small $$\textstyle\left( \begin{array}{c} \delta_1 \\
\delta_2 \\ \delta_3
\end{array} \right) =
\left( \begin{array}{ccc}
 x^2+\frac{5}{4} xy & \frac{3}{4} xy+y^2 &\frac{1}{4} xz^2-\frac{1}{4}xz \\
 0 & 0  & (xz+y) \\
 4xy^2+y^3+25x^2& 3y^3-x^2+20xy &y^2z^2-5xz-y^2z+x
\end{array} \right)
\left( \begin{array}{c} \dx \\ \dy \\
\dz
\end{array} \right)$$ }

\noindent On v\'erifie que:\smallskip

\noindent -) $\delta_1(h)= (\frac{1}{4}xz+\frac{19}{4}x+6y)h,\quad
\delta_2(h)=xh, \quad \delta_3(h)= (y^2z+19y^2+120x)h$,\\
-) le d\'eterminant de la matrice  des coefficients pr\'ec\'edente
vaut $h$,\\
-) $ [\delta_1,\delta_2]=(x+y-\frac{1}{4}xz)\delta_2,\quad
 \quad
[\delta_2,\delta_3]=(y^2z-4y^2-25x)\delta_2$,\\
-) $
 [\delta_1,\delta_3]=-(3y^2+30x)\delta_1
 +(\frac{5}{2}xz-\frac{1}{2}x)\delta_2+(\frac{5}{4}x+\frac{7}{4}y)\delta_3$.
\medskip

Nous allons prouver que $D$ n'est pas un diviseur de Spencer (au
voisinage du point $p=(0,0,0))$. Pour cela voyons que la fibre en
$p$ du complexe
$$ \DX\Lotimes_{\VCERO} \OX = \DX\otimes_{\VCERO}\SP_{\VCERO}$$
a une cohomologie non nulle en degr\'e $-1$.
\medskip

Notons encore $d^\bullet$ la diff\'erentielle du complexe
$\DX\otimes_{\VCERO}\SP_{\VCERO}$. L'image de $d^{-2}$ est
engendr\'ee comme $\DX$-module \`a gauche par:\smallskip

\noindent $ d^{-2}(1\otimes(\delta_1\wedge\delta_2))  =
-\delta_2\otimes\delta_1 +
(\delta_1-x-y+\frac{1}{4}xz)\otimes\delta_2$,\smallskip

\noindent $d^{-2}(1\otimes(\delta_1\wedge\delta_3) =
(-\delta_3+3y^2+30x)\otimes\delta_1+(-\frac{5}{2}xz+\frac{1}{2}x)\otimes\delta_2+
 (\delta_1-\frac{5}{4}x-\frac{7}{4}y)\otimes\delta_3$,\smallskip

\noindent $d^{-2}(1\otimes(\delta_2\wedge\delta_3) =
(-\delta_3-y^2z+4y^2+25x)\otimes\delta_2+\delta_2\otimes\delta_3$.
\medskip

Par cons\'equent, tout \'el\'ement $P_1\otimes\delta_1 +
P_2\otimes\delta_2+P_3\otimes \delta_3$ de l'image de $d^{-2}$
v\'erifie que $P_3$ appartient \`a l'id\'eal \`a gauche $I$ de $\DX$
engendr\'e par $\delta_2$ et
$\zeta=\delta_1-\frac{5}{4}x-\frac{7}{4}y$.
\medskip

Nous allons prouver que l'image de $d^{-2}$ est contenue
strictement dans le noyau de $d^{-1}$ au point $p=(0,0,0)$.
\medskip

  Consid\'erons les op\'erateurs diff\'erentiels:
   {\small
 $$Q_1=-{y} z^{3} {\dz}^{2}-y^{2} {z} {\dx} {\dz}-3
y^{2} {z} {\dy} {\dz}+{y} z^{2} {\dz}^{2}-2 {y} z^{2} {\dz}+4
y^{2} {\dx} {\dz}-25 {x} {z} {\dy} {\dz}+8 {y} {z} {\dz}+$$
 $$25 {x} {\dx} {\dz}-{x} {\dy} {\dz}-5 {y} {\dy} {\dz}-5 {z}
{\dz}^{2}+{\dz}^{2}-60 {\dz}, $$

$$Q_2=\frac{1}{4}\left(yz^4\dz^2+4xyz^2\dx\dz+6y^2z^2\dx\dz+{2}y^2z^2\dy\dz-{2}yz^3\dz^2+4xy^2\dx^2+{5}y^3\dx^2-\right.$$
$$4xy^2\dx\dy-{2}y^3\dx\dy-{3}y^3\dy^2+4yz^3\dz-4xyz\dx\dz-6y^2z\dx\dz-{2}y^2z\dy\dz+$$
$${25}xz^2\dy\dz+yz^2\dz^2+4xyz\dx+6y^2z\dx+{2}y^2z\dy+{25}xy\dx\dy+4x^2\dy^2-xy\dy^2+$$
$$20y^2\dy^2-24yz^2\dz-{5}xz\dy\dz+20yz\dy\dz+{2}yz^2-60xy\dx-16y^2\dx-44y^2\dy+{25}xz\dy+$$
$$\left.20yz\dz-4x\dy\dz-4y\dy\dz-20yz+400x\dx-16x\dy+340 y\dy-{45}z\dz+60y+{9}\dz-365\right),$$

$$ Q_3=\frac{1}{4}\left(4
{x} {z} {\dy} {\dz}+4 {y} {z} {\dy} {\dz}-z^{2} {\dz}^{2}-4 {x} {\dx}
{\dz}-5 {y} {\dx} {\dz}+{y} {\dy} {\dz}+{z} {\dz}^{2}-13 {z} {\dz}+10
{\dz}\right).$$ }
\medskip

On v\'erifie que $Q_1\delta_1+Q_2\delta_2+Q_3\delta_3=0$ et donc
$Q=Q_1\otimes\delta_1 + Q_2\otimes\delta_2+Q_3\otimes \delta_3$
est dans le noyau de $d^{-1}$.
\medskip

Voyons que le germe en $p$ de $Q_3$ n'appartient pas \`a $I_p$, et
donc le germe en $p$ de $Q$ n'appartient pas \`a la fibre en $p$
de l'image de $d^{-2}$.
\medskip

Notons $R =\C[x,y,z]\subset \OO_{X,p}$, ${\mathfrak m} =
R(x,y,z)$, $S = R[\xi_1,\xi_2,\xi_3]$, l'extension $S_{\mathfrak
m} = R_{\mathfrak m}[\xi_1,\xi_2,\xi_3] \subset\Gr \D_{X,p}=
\OO_{X,p}[\xi_1,\xi_2,\xi_3]$ \'etant fid\`element plate.
\medskip

On a $ [\delta_2,\zeta]= (x+y-\frac{1}{4}xz)\delta_2$, d'o\`u le
$\OX$-module engendr\'e par $\delta_2$ et $\zeta$ est un
sous-alg\'ebro\"{\i}de de Lie de $F^1\DX$. D'autre part, les symboles
$\sigma(\delta_2)$ et $\sigma(\zeta)=\sigma(\delta_1)$ forment une
suite r\'eguli\`ere dans $R$, et par platitude dans $\Gr \D_{X,p}$.
\medskip

On v\'erifie que $(S(\sigma(\zeta),\sigma(\delta_2)):\sigma(Q_3))
= S(x,y)$, d'o\`u
$(S(\sigma(\zeta),\sigma(\delta_2)):\sigma(Q_3))\cap R = R(x,y)$,
et donc $\sigma(Q_3)\notin R_{\mathfrak m}[\xi_1,\xi_2,\xi_3]
(\sigma(\zeta),\sigma(\delta_2))$.
 Par fid\`ele platitude, $\sigma (Q_3) \notin \Gr\D_{X,p} (\sigma(\zeta),\sigma(\delta_2))$.
\medskip

Le r\'esultat cherch\'e est donc une cons\'equence de la
proposition suivante, dont la preuve est la m\^eme que celle de la
proposition 4.1.2 de \cite{calde_ens}.
\end{ejemplo}

\begin{proposicion} Soit $I$ un id\'eal \`a gauche de $\D_{X,p}$ engendr\'e par des op\'erateurs
$P_1,\ldots,P_t$ d'ordre $\leq 1$ tels que:
\begin{itemize}
 \item La suite des symboles $\{\sigma(P_1),\ldots,\sigma(P_t)\}$
 est r\'eguli\`ere dans $\Gr \D_{X,p}$.
 \item
 $ [P_i,P_j]\in\sum_{k=1}^t\OX P_k,\ 1\leq i,j\leq t$.
\end{itemize}
Alors $\sigma(I)=\Gr \D_{X,p}(\sigma(P_1),\ldots,\sigma(P_t))$.
\end{proposicion}

\begin{ejemplo}
Dans le cas non localement quasi-homog\`ene il y a des exemples de
connexions logarithmiques int\'egrables telles que le morphisme
naturel
\begin{equation}
\label{eq:comp-cli}\DX\Lotimes_{\VCERO} \E(kD)\to \E(\star
D)\end{equation} n'est un isomorphisme pour aucun $k> 0$. Prendre par
exemple une courbe plane $D\subset X=\C^2$, d'\'equation $f=0$, qui
n'est pas quasi-homog\`ene et $\E = \OX$. Dans ce cas on a:
\begin{enumerate}
\item[1)] Les complexes $\DX\Lotimes_{\VCERO} \OX(kD)$, $k\in \ZZ$,  sont
concentr\'es en degr\'e $0$, car $D$ est de Koszul,
\item[2)] les
morphismes $ \DX\otimes_{\VCERO} \OX(kD)\to \OX(\star D)$, $k>0$,
sont surjectifs, car les fonctions m\'eromorphes sont engendr\'ees
sur $\DX$ par $f^{-1}$ dans le cas des courbes planes,
\item[3)] le morphisme $\DX\Lotimes_{\VCERO} \OX(D)\to \OX(\star D)$
n'est pas un isomorphisme, d'apr\`es \cite{calde_mond_nar_cas} et le
th\'eor\`eme \ref{teor:gordo},
\item[4)] tous les $\DX$-modules $\DX\otimes_{\VCERO} \OX(kD)$, $k\in
\ZZ$, ont m\^eme cycle caract\'eristique, car $D$ est de Koszul (voir
la preuve du th\'eor\`eme \ref{tcl}).
\end{enumerate}
Par cons\'equent, le morphisme $\DX\Lotimes_{\VCERO} \OX(kD)\to
\OX(\star D)$ n'est un isomorphisme pour aucun $k>0$.

Pourtant, dans le cas localement quasi-homog\`ene il est raisonnable
de se demander si le morphisme (\ref{eq:comp-cli}) est un
isomorphisme pour $k\gg 0$. Dans un travail en cours nous \'etudions
cette question en partant de la g\'en\'eralisation des r\'esultats de
\cite{calde_nar_compo} au cas des connexions logarithmiques
int\'egrables.
\end{ejemplo}

\begin{problema} Est-ce que le complexe de de Rham logarithmique
$\OmX(\log D)$ de tout diviseur libre $D\subset X$ est analytiquement
constructible?, ou de fa\c{c}on \'equivalente d'apr\`es le th\'eor\`eme
de constructibilit\'e de Kashiwara et \cite{meb_formalisme}, chap.
II, th. (4.1.5), est-ce que le
complexe $\DX\Lotimes_{\VCERO} \OX$ est \`a cohomologie holonome?\\
On peut aussi se poser la m\^eme question pour les complexes de de
Rham logarithmiques associ\'es \`a des connexions logarithmiques
int\'egrables arbitraires.
\end{problema}

\begin{problema} Est-ce que si $\OX$ est addmissible (i.e. $D$ est un diviseur de Spencer),
alors toute
connexion logarithmique int\'egrable est admissible?
\end{problema}

\begin{problema} Soit $D$ un diviseur libre qui satisfait le th\'eor\`eme de comparaison
logarithmique. Dans \cite{calde_mond_nar_cas} on a conjectur\'e que
$D$ est Euler-homog\`ene et d'apr\`es le corollaire \ref{criterio}
nous savons que $D$ est de Spencer. Pourtant, il existe des diviseurs
libres Euler-homog\`enes de Spencer, et m\^eme de Koszul, qui ne
satisfont pas le th\'eor\`eme de comparaison logarithmique: il suffit
de prendre une surface d'\'equation $zf(x,y)=0$ o\`u $f(x,y)=0$ est
l'\'equation r\'eduite d'une courbe plane non quasi-homog\`ene. Une
question naturelle est de trouver des conditions suffisantes sur les
diviseurs libres Euler-homog\`enes, autres que la
quasi-homog\'en\'eit\'e locale (resp. que d'\^etre de Koszul) qui
garantissent le th\'eor\`eme de comparaison logarithmique (resp.
d'\^etre de Spencer). \end{problema}

\section*{Appendice} \renewcommand{\thesection}{A} \setcounter{numero}{0}

Dans cet appendice on fixe un homomorphisme d'anneaux commutatifs $k
\to A$ (i.e., $A$ est une $k$-alg\`ebre). \'Etant donn\'ee une
$(k,A)$-alg\`ebre de Lie $L$, notons $\U(L)$ son alg\`ebre
enveloppante cf. \cite{rine-63}. Rappelons qu'elle est construite
comme le quotient de l'alg\`ebre tensorielle $\Tens(A\oplus L)$ par
l'id\'eal bilat\`ere engendr\'e par les \'el\'ements:\medskip

\noindent $c\cdot 1 - i(c), c\in k$,\\
$i(a)\otimes i(b) - i(ab), a,b\in A\subset A\oplus L$,\\
$i(\lambda)\otimes i(a) - i(a)\otimes i(\lambda) - i(\lambda(a)),
\lambda\in L, a\in A$,\\
$i(\lambda)\otimes i(\mu) - i(\mu)\otimes i(\lambda) -
i([\lambda,\mu]), \lambda, \mu \in L$,\\
$i(a\lambda) - i(a)\otimes i(\lambda), \lambda\in L, a\in A$,\medskip

\noindent $i: A\oplus L \to \Tens(A\oplus L)$ \'etant l'inclusion
naturelle.
\medskip

Une description alternative est donn\'ee dans \cite{calde_ens}, prop.
2.2.5, en utilisant l'alg\`ebre tensorielle des bimodules et le fait
que $A\oplus L$ est muni d'une structure de $(A,A)$-bimodule:
$$ a(b+\lambda) = ab + a\lambda,\quad (b+\lambda)a= [ba + \lambda(a)] +
a\lambda,\quad a,b\in A, \lambda\in L.$$De plus, $A\oplus L$ a aussi
une structure de $(k,A)$-alg\`ebre de Lie.
\medskip

L'alg\`ebre $\U(L)$ est un anneau filtr\'e de la mani\`ere
\'evidente.
\medskip

Si $L$ est un module libre de rang $n$ sur $A$, le th\'eor\`eme de
Poincar\'e-Birkhoff-Witt (\cite{rine-63}, th. 3.1) nous dit que le
gradu\'e associ\'e \`a l'anneau $\U(L)$ est canoniquement isomorphe
\`a l'alg\`ebre sym\'etrique du $A$-module $L$.
\medskip

Consid\'erons maintenant un morphisme de $(k,A)$-alg\`ebres de Lie
$L_0 \to L$ et le morphisme induit sur les alg\`ebres enveloppantes
correspondantes $U_0= \U(L_0)\to \U(L)=U$. Si $\varepsilon\in L_0$,
on \'ecrira aussi $\varepsilon$ pour son image dans $L$.

Supposons que $M$ (resp. $N$) est un $U$-module \`a droite (resp. un
$U_0$-module \`a gauche). Alors $U\otimes_{U_0} N$ est un $U$-module
\`a gauche, et donc $M\otimes_A[U\otimes_{U_0} N]$ est un $U$-module
\`a droite.

D'autre part, comme $M$ est aussi un $U_0$-module \`a droite par
restriction des scalaires, $M\otimes_A N$ est un $U_0$-module \`a
droite.

Il est clair que le morphisme:
$$ \sigma: M\otimes_A N \to M\otimes_A[U\otimes_{U_0} N],\quad
\sigma(m\otimes n) = m\otimes[1\otimes n]$$est $U_0$-lin\'eaire \`a
droite, et induit donc un morphisme $U$-lin\'eaire \`a droite:
$$\lambda: [M\otimes_A N]\otimes_{U_0} U \to M\otimes_A[U\otimes_{U_0}
N].$$
\medskip

\begin{teorema} \label{A.1}
Pour tout $U$-module \`a droite $P$ et pour tout morphisme
$U_0$-lin\'eaire \`a droite $\alpha : M\otimes_A N \to P$, il existe
un seul morphisme $U$-lin\'eaire \`a droite $\beta:
M\otimes_A[U\otimes_{U_0} N] \to P$ tel que $\beta \circ \sigma =
\alpha$.
\end{teorema}
\medskip

Notons que le th\'eor\`eme \ref{A.1} nous dit que le morphisme
$\sigma$ v\'erifie la m\^eme propri\'et\'e universelle que le
morphisme naturel $[M\otimes_A N]\to [M\otimes_A N]\otimes_{U_0} U$.
On a donc le corollaire suivant:
\medskip

\begin{corolario} \label{A.2} Sous les hypoth\`eses pr\'ec\'edentes, le morphisme
$$\lambda: [M\otimes_A N]\otimes_{U_0} U \to M\otimes_A[U\otimes_{U_0}
N]$$ est un isomorphisme de $U$-modules \`a droite.
\end{corolario}
\medskip

\noindent {\bf Preuve du th\'eor\`eme \ref{A.1}:}  L'unicit\'e est
claire.
\medskip

Notons $\Delta_0 = A\oplus L_0,  \Delta = A\oplus L$, $T_0 =
\Tens(\Delta_0), T = \Tens(\Delta)$.

Le morphisme $\beta$ proviendra d'une application $k$-multilin\'eaire
$$\nu: M \times T \times N \to P$$satisfaisant certaines conditions.

Comme $T = \bigoplus T^r$, avec $T^0 = k$ et $T^r = \Delta^{\otimes
r}$ (produit tensoriel sur $k$), il suffit de d\'efinir des
applications $k$-multilin\'eaires $$ \nu_r : M\times T^r \times N \to
P, r\geq 0.$$ Pour $r=0$, on doit avoir clairement $$\nu_0(m,c,n) =
\alpha((mc)\otimes n).$$ On d\'efinit de mani\`ere r\'ecurrente
$\nu_r:M\times T^r \times N \to P$,
$$\nu_r(m,t_1\otimes \cdots \otimes t_r,n) =
\nu_{r-1}(mt_1,t_2\otimes \cdots \otimes t_r,n) -
\nu_{r-1}(m,t_2\otimes \cdots \otimes t_r,n)\lambda_1,$$ o\`u on a
\'ecrit $t_i = a_i + \lambda_i$, $a_i\in A, \lambda_i\in L$.

On voit facilement que l'expression pr\'ec\'edente est bien d\'efinie
sur le produit tensoriel et que:
$$\nu_r(ma,t_1\otimes \cdots \otimes t_r,n) = \nu_r(m,(at_1)\otimes
\cdots \otimes t_r,n), a\in A.$$ Notons $I_r =\{1,\dots, r\}$, et si
$E\subset I_r$, $E=\{i_1 < \cdots < i_l\}$ on \'ecrira: $$ t_E =
t_{i_1} \cdots t_{i_l},\quad \lambda_E = \lambda_{i_l}
\lambda_{i_{l-1}} \cdots \lambda_{i_1}\quad\text{(attention \`a
l'ordre)}.$$ Si $a\in A$ on notera $\lambda_E(a) = \lambda_{i_l}(
\lambda_{i_{l-1}}( \cdots (\lambda_{i_1}(a))\cdots ))$.
\medskip

\noindent On d\'emontre par r\'ecurrence la formule suivante: $$
\nu_r(m,t_1\otimes \cdots \otimes t_r,n) = \sum_{E\subset I_r}
(-1)^{\sharp (I_r - E)} \alpha((m t_E)\otimes n)\lambda_{I_r-E}.$$
\`A partir du lemme de commutation \ref{A.3}, on d\'emontre les
relations suivantes:
$$ \nu_r(ma,t_1\otimes \cdots \otimes t_r,n) = \nu_r(m,t_1\otimes
\cdots \otimes t_r,n)a, a\in A,$$ $$ \nu_r(m,t_1\otimes \cdots
\otimes (t_ra),n) = \nu_r(m,t_1\otimes \cdots \otimes t_r,an), a\in
A.$$

La preuve des trois lemmes suivants est laiss\'ee au lecteur.
\smallskip

\begin{lema} \label{A.3} {\em (Lemme de commutation)}  Pour $t_i = a_i + \lambda_i$, $i\in
I_r$, $a\in A$ et $E\subset I_r$ on a: $$ a t_E = \sum_{E'\subset E}
(-1)^{\sharp(E-E')} t_{E'} \lambda_{E-E'}(a),$$ $$ \lambda_E a =
\sum_{E'\subset E} \lambda_{E'}(a) \lambda_{E-E'}.$$
\end{lema}\medskip

\begin{lema} \label{A.4-app}
 Pour $m\in M, n\in N, t_i\in \Delta$, $i=1,\dots,r-1$, $t_r\in
\Delta_0$, on a: $$ \nu_r(m,t_1\otimes \cdots \otimes t_r,n) =
\nu_{r-1}(m,t_1\otimes \cdots \otimes t_{r-1},t_r n).$$
\end{lema}\medskip

\begin{lema} \label{A.5}
Pour $m\in M, n\in N, t_i\in \Delta$, $i=2,\dots,r$, $\lambda_1\in
L$, on a: $$ \nu_{r-1}(m\lambda_1,t_2\otimes \cdots \otimes t_r,n) -
\nu_{r}(m,\lambda_1\otimes t_2\otimes \cdots \otimes t_{r},n)=$$ $$ =
\nu_{r-1}(m,t_2\otimes \cdots \otimes t_r,n)\lambda_1.$$
\end{lema}\medskip

Ensuite on voit que l'application $\nu$ passe au quotient $U = T/J$
et d\'efinit une nouvelle application $k$-multilin\'eaire not\'ee
encore
$$\nu: M \times U \times N \to P,$$ qui gr\^ace aux lemmes pr\'ec\'edents
d\'efinit l'application $U$-lin\'eaire \`a droite cherch\'ee $$
\beta: M\otimes_A[U\otimes_{U_0} N] \to P$$et le th\'eor\`eme
\ref{A.1} est d\'emontr\'e.\hfill $\Box$ 
\medskip

Soit maintenant $M$ (resp. $N$) un $U$-module \`a gauche (resp. un
$U_0$-module \`a gauche). Alors $U\otimes_{U_0} N$ est un $U$-module
\`a gauche, et donc $[U\otimes_{U_0} N]\otimes_A M$ est un $U$-module
\`a gauche.

D'autre part, comme $M$ est aussi un $U_0$-module \`a gauche par
restriction des scalaires, $N\otimes_A M$ est un $U_0$-module \`a
gauche. Il est clair que le morphisme
$$ n\otimes m \in N\otimes_A M \mapsto [1\otimes n]\otimes m \in [U\otimes_{U_0} N]\otimes_A
M$$est $U_0$-lin\'eaire \`a gauche et induit donc un morphisme
$U$-lin\'eaire \`a gauche:
$$ \lambda': U\otimes_{U_0} [N\otimes_A M] \xrightarrow{} [U\otimes_{U_0} N]\otimes_A
M.$$ Le th\'eor\`eme suivant se d\'emontre de fa\c{c}on tout \`a fait
analogue au th\'eor\`eme \ref{A.1} et au corollaire \ref{A.2}:
\medskip

\begin{teorema}\label{A.6} Sous les hypoth\`eses pr\'ec\'edentes, le
morphisme
$$ \lambda': U\otimes_{U_0} [N\otimes_A M] \xrightarrow{} [U\otimes_{U_0} N]\otimes_A
M$$
 est un isomorphisme de $U$-modules \`a gauche.
\end{teorema}

La situation pr\'ec\'edente s'applique au cas o\`u $L_0 = 0$ et $U_0
= A \to U$ est l'inclusion:

\begin{corolario} \label{A.7} Soit $M$ un $U$-module \`a gauche et $N$
un $A$-module. Alors il existe un unique isomorphisme $U$-lin\'eaire
\`a gauche
$$ \lambda': U\otimes_{A} [N\otimes_A M] \xrightarrow{} [U\otimes_{A} N]\otimes_A
M$$qui envoie $1\otimes[n\otimes m]$ dans $[1\otimes n]\otimes m$.
\end{corolario}

\def\cprime{$'$}

\bigskip

{\small \noindent Departamento de \'Algebra,
 Facultad de  Matem\'aticas, Universidad de Sevilla,\\ PO Box 1160, 41080
 Sevilla, Spain}. \\
{\small {\it E-mail}:  $\{$calderon,narvaez$\}$@algebra.us.es
 }

\end{document}